\documentclass[10pt]{article}
\usepackage{mathrsfs}
\usepackage{amsfonts}
\usepackage{amsmath}
\usepackage{amsfonts,amssymb}
\usepackage{mathrsfs}
\usepackage{amssymb}
\allowdisplaybreaks

\numberwithin{equation}{section}

\date{}

\begin{document}
\title{A result on fractional $(a,b,k)$-critical covered graphs
}
\author{\small  Sizhong Zhou\footnote{Corresponding author. E-mail address: zsz\_cumt@163.com (S. Zhou)}, Quanru Pan\\
\small  School of Science, Jiangsu University of Science and Technology\\
\small  Mengxi Road 2, Zhenjiang, Jiangsu 212003, P. R. China\\
}

\maketitle
\begin{abstract}
\noindent For a graph $G$, the set of vertices in $G$ is denoted by $V(G)$, and the set of edges in $G$ is denoted by $E(G)$.
A fractional $[a,b]$-factor of a graph $G$ is a function $h$ from $E(G)$ to $[0,1]$ satisfying $a\leq d_G^{h}(v)\leq b$ for
every vertex $v$ of $G$, where $d_G^{h}(v)=\sum\limits_{e\in E(v)}{h(e)}$ and $E(v)=\{e=uv:u\in V(G)\}$. A graph $G$ is called
fractional $[a,b]$-covered if $G$ contains a fractional $[a,b]$-factor $h$ with $h(e)=1$ for any edge $e$ of $G$. A graph $G$
is called fractional $(a,b,k)$-critical covered if $G-Q$ is fractional $[a,b]$-covered for any $Q\subseteq V(G)$ with $|Q|=k$.
In this article, we demonstrate a neighborhood condition for a graph to be fractional $(a,b,k)$-critical covered. Furthermore,
we claim that the result is sharp.
\\
\begin{flushleft}
{\em Keywords:} graph; neighborhood; fractional $[a,b]$-factor; fractional $[a,b]$-covered graph; fractional $(a,b,k)$-critical
covered graph.

(2010) Mathematics Subject Classification: 05C70
\end{flushleft}
\end{abstract}

\section{Introduction}

We may use a graph to model a real-world network by considering the nodes of the network as the vertices of the graph and
considering the links between the nodes in the network as the edges of the graph. For example, in a communication network
we use nodes to represent cities and use links to act for channels, or in the World Wide Web we apply nodes to stand for
Web pages and apply links to represent hyperlinks between Web pages. Henceforth, we use the term ``graph" instead of ``network".

The fractional factor problems in graphs can be considered as a relaxation of the well-known cardinality matching problems.
The fractional factor problem admits widespread applications in many areas such as scheduling, network design, and combinatorial
polyhedron because such real situations can be modelled using dynamical systems defined by graphs, see for more information \cite{AMV}.
For instance, in a communication network, several large data packets can be sent to various destinations through several channels.
If we allow the large data packets to be partitioned into small parcels, then the efficiency of the network can be improved. The
feasible assignment of data packets can be seen as a fractional flow problem, and it becomes a fractional factor problem while the
destinations and the sources of a network are disjoint \cite{LYZ0}.

In the graph theory model, we delete the vertices of a graph corresponding to the damaged sites of a network and occupy the edges
of a graph corresponding to the links of a network, and we look for the corresponding fractional factor in the resulting graph. Hence,
the data transmission problem at the moment in a communication network can be converted into the related problem on the existence
of the fractional factor critical graph, the fractional factor covered graph or the fractional factor critical covered graph in the
corresponding graph of the network. The neighborhood condition of a graph is often used to measure the robustness and vulnerability
of a network, which is an important parameter in data transmission and network design, and so on.

All graphs discussed are assumed to be finite, undirected and simple. For a graph $G$, the set of vertices in $G$ is
denoted by $V(G)$, and the set of  edges in $G$ is denoted by $E(G)$. For $v\in V(G)$, the degree of $v$ in $G$ is
the number of edges of $G$ incident with $v$ and is denoted by $d_G(v)$, and the neighborhood of $v$ in $G$ is the
set of vertices of $G$ adjacent to $v$ and is denoted by $N_G(v)$. Distinctly, $d_G(v)=|N_G(v)|$. For any subset
$X$ of $V(G)$, the subgraph of $G$ induced by $X$ is denoted by $G[X]$ and the neighborhood of $X$ in $G$ is denoted
by $N_G(X)$. Setting $G-X=G[V(G)\setminus X]$. A subset $X$ of $V(G)$ is called independent if $G[X]$ does not possess
edges. The minimum degree of $G$ is defined by $\delta(G)=\min\{d_G(v):v\in V(G)\}$. Assume that $c$ is a real number.
Recalling that $\lfloor c\rfloor$ is the greatest integer satisfying $\lfloor c\rfloor\leq c$.

For positive integers $a$ and $b$ satisfying $a\leq b$, an $[a,b]$-factor of a graph $G$ is a spanning subgraph $F$ of
$G$ satisfying $a\leq d_F(v)\leq b$ for all $v\in V(G)$. Assume that $a=b=r$, then we call it an $r$-factor. A fractional
$[a,b]$-factor of a graph $G$ is a function $h$ from $E(G)$ to $[0,1]$ satisfying $a\leq d_G^{h}(v)\leq b$ for every
vertex $v$ of $G$, where $d_G^{h}(v)=\sum\limits_{e\in E(v)}{h(e)}$ and $E(v)=\{e=uv:u\in V(G)\}$. If $a=b=r$, then we
call it a fractional $r$-factor. A graph $G$ is called fractional $[a,b]$-covered if $G$ contains a fractional
$[a,b]$-factor $h$ with $h(e)=1$ for any edge $e$ of $G$. If $a=b=r$, then we call $G$ being a fractional $r$-covered
graph. A graph $G$ is called fractional $(a,b,k)$-critical covered if $G-Q$ is fractional $[a,b]$-covered for any
$Q\subseteq V(G)$ with $|Q|=k$, where $k$ is a nonnegative integer. When $a=b=r$, the fractional $(a,b,k)$-critical
covered graph is the fractional $(r,k)$-critical covered graph.

The previous works have implied that there is a close relationship between neighborhood and the existence of factors
and fractional factors in graphs. Amahashi and Kano \cite{AK} and Berge and Las Vergnas \cite{BV} derived a neighborhood
condition for graphs having $[1,b]$-factors, independently. Kano \cite{K} posed a neighborhood condition for the
existence of $[a,b]$-factors in graphs. Zhou, Pu and Xu \cite{ZPX} claimed a sharp neighborhood condition for a graph
possessing a fractional $r$-factor. More results on factors and fractional factors of graphs were gained by Xiong \cite{X},
Axenovich and Rollin \cite{AR}, Zhou and Sun \cite{ZS}, Zhou, Sun and Xu \cite{ZSX}, Sun and Zhou\cite{SZ}, Zhou, Sun and Ye
\cite{ZSY}, Zhou \cite{Zr,Zhr,Zhs,Zs}, Zhou, Xu and Sun \cite{ZXS}, Zhou, Yang and Xu \cite{ZYX}, Zhou, Zhang and Xu \cite{ZTX},
Gao, Guirao and Wu \cite{GGW}, Gao and Wang \cite{GWa,GWn}, Yuan and Hao \cite{YH}, Chiba and Yamashita \cite{CY}.

A neighborhood condition for the existence of fractional $r$-factors in graphs was demonstrated by Zhou, Pu and Xu
\cite{ZPX}, which is the following result.

\medskip

\noindent{\textbf{Theorem 1}} (\cite{ZPX}). Let $r\geq1$ be an integer, and let $G$ be a graph of order $n$ with
$n\geq6r-12+\frac{6}{r}$. Assume, for every subset $X$ of $V(G)$, that
$$
N_G(X)=V(G) \ \ \ \ \ \ if \ \ \ \ \ \ |X|\geq\bigg\lfloor\frac{rn}{2r-1}\bigg\rfloor; \ \ \ or
$$
$$
|N_G(X)|\geq\frac{2r-1}{r}|X| \ \ \ \ \ \ if \ \ \ \ \ \ |X|<\bigg\lfloor\frac{rn}{2r-1}\bigg\rfloor.
$$
Then $G$ possesses a fractional $r$-factor.

\medskip

In this article, we generalize Theorem 1, and claim a neighborhood condition for graphs being fractional
$(a,b,k)$-critical covered graphs.

\medskip

\noindent{\textbf{Theorem 2.}} Let $a,b$ and $k$ be integers with $b\geq a\geq2$ and $k\geq0$, and let $G$ be a graph of
order $n$ with $n\geq\frac{(a+b-2)(2a+b-3)+2}{b}+\frac{bk}{b-1}$. Assume, for every vertex subset $X$ of $G$, that
$$
N_G(X)=V(G) \ \ \ \ \ \ if \ \ \ \ \ \ |X|\geq\bigg\lfloor\frac{(b(n-1)-bk)n-2(n-1)}{(a+b-1)(n-1)}\bigg\rfloor; \ \ \ or
$$
$$
|N_G(X)|\geq\frac{(a+b-1)(n-1)}{b(n-1)-bk-2}|X| \ \ \ \ \ \ if \ \ \ \ \ \ |X|<\bigg\lfloor\frac{(b(n-1)-bk)n-2(n-1)}{(a+b-1)(n-1)}\bigg\rfloor.
$$
Then $G$ is fractional $(a,b,k)$-critical covered.

\medskip

The following result holds for $k=0$ in Theorem 2.

\medskip

\noindent{\textbf{Corollary 1.}} Let $a,b$ be integers with $b\geq a\geq2$, and let $G$ be a graph of order $n$ with
$n\geq\frac{(a+b-2)(2a+b-3)+2}{b}$. Assume, for every vertex subset $X$ of $G$, that
$$
N_G(X)=V(G) \ \ \ \ \ \ if \ \ \ \ \ \ |X|\geq\bigg\lfloor\frac{bn-2}{a+b-1}\bigg\rfloor; \ \ \ or
$$
$$
|N_G(X)|\geq \frac{(a+b-1)(n-1)}{b(n-1)-2}|X| \ \ \ \ \ \ if \ \ \ \ \ \ |X|<\bigg\lfloor\frac{bn-2}{a+b-1}\bigg\rfloor.
$$
Then $G$ is fractional $[a,b]$-covered.

\medskip

We can acquire the following result when $a=b=r$ in Theorem 2.

\medskip

\noindent{\textbf{Corollary 2.}} Let $r$ and $k$ be integers with $r\geq2$ and $k\geq0$, and let $G$ be a graph of
order $n$ with $n\geq6r-12+\frac{8}{r}+\frac{rk}{r-1}$. Assume, for every vertex subset $X$ of $G$, that
$$
N_G(X)=V(G) \ \ \ \ \ \ if \ \ \ \ \ \ |X|\geq\bigg\lfloor\frac{(r(n-1)-rk)n-2(n-1)}{(2r-1)(n-1)}\bigg\rfloor; \ \ \ or
$$
$$
|N_G(X)|\geq \frac{(2r-1)(n-1)}{r(n-1)-rk-2}|X| \ \ \ \ \ \ if \ \ \ \ \ \ |X|<\bigg\lfloor\frac{(r(n-1)-rk)n-2(n-1)}{(2r-1)(n-1)}\bigg\rfloor.
$$
Then $G$ is fractional $(r,k)$-critical covered.

\medskip

The following result is obtained if $a=b=r$ in Corollary 1.

\medskip

\noindent{\textbf{Corollary 3.}} Let $r$ be an integer with $r\geq2$, and let $G$ be a graph of order $n$ with
$n\geq6r-12+\frac{8}{r}$. Assume, for every vertex subset $X$ of $G$, that
$$
N_G(X)=V(G) \ \ \ \ \ \ if \ \ \ \ \ \ |X|\geq\bigg\lfloor\frac{rn-2}{2r-1}\bigg\rfloor; \ \ \ or
$$
$$
|N_G(X)|\geq \frac{(2r-1)(n-1)}{r(n-1)-2}|X| \ \ \ \ \ \ if \ \ \ \ \ \ |X|<\bigg\lfloor\frac{rn-2}{2r-1}\bigg\rfloor.
$$
Then $G$ is fractional $r$-covered.

\medskip

\section{The proof of Theorem 2}

The proof of Theorem 2 relies on the following theorem, which is a special case of fractional $(g,f)$-covered graph
theorem obtained by Li, Yan and Zhang \cite{LYZ}.

\medskip

\noindent{\textbf{Theorem 3}} (\cite{LYZ}). Let $a,b$ be nonnegative integers satisfying $b\geq a$, and let $G$ be a graph.
Then $G$ is fractional $[a,b]$-covered if and only if
$$
\theta_G(S,T)=b|S|+d_{G-S}(T)-a|T|\geq\varepsilon(S)
$$
for any vertex subset $S$ of $G$ and $T=\{x:x\in V(G)\setminus S, d_{G-S}(x)\leq a\}$, where $\varepsilon(S)$ is defined by
\[
 \varepsilon(S)=\left\{
\begin{array}{ll}
2,&if \ S \ is \ not \ independent,\\
1,&if \ S \ is \ independent, \ and \ there \ is \ an \ edge \ joining \ S \ and \ V(G)\setminus(S\cup T) \ or\\
&there \ is \ an \ edge \ e=uv \ joining \ S \ and \ T \ satisfying \ d_{G-S}(v)=a \ for \ v\in T,\\
0,&otherwise.\\
\end{array}
\right.
\]

\medskip

\noindent{\it Proof of Theorem 2.} \ Let $D\subseteq V(G)$ with $|D|=k$, and let $H=G-D$. It suffices to claim that $H$
is fractional $[a,b]$-covered. Assume that $H$ is not fractional $[a,b]$-covered. Then using Theorem 3, we possess
$$
\theta_H(S,T)=b|S|+d_{H-S}(T)-a|T|\leq\varepsilon(S)-1,\eqno(1)
$$
for a vertex subset $S$ of $H$, where $T=\{x:x\in V(H)\setminus S, d_{H-S}(x)\leq a\}$. Evidently, $T\neq\emptyset$ by (1)
and $\varepsilon(S)\leq|S|$. Thus, we define
$$
\beta=\min\{d_{H-S}(x):x\in T\}.
$$
By virtue of the definition of $T$, we acquire
$$
0\leq\beta\leq a.
$$

\noindent{\it Claim 1.} \ $\delta(H)\geq\frac{(a-1)n+b-(a-1)k+2}{a+b-1}$.

\noindent{\it Proof.} \ Let $v\in V(G)$ such that $d_G(v)=\delta(G)$, and set $W=V(G)\setminus N_G(v)$. Distinctly,
$v\notin N_G(W)$, and so $N_G(W)\neq V(G)$. Then using the condition of Theorem 2, we deduce
$$
|N_G(W)|\geq\frac{(a+b-1)(n-1)}{b(n-1)-bk-2}|W|.
$$
Note that $|W|=n-\delta(G)$ and $|N_G(W)|\leq n-1$. Thus, we derive
$$
n-1\geq|N_G(W)|\geq\frac{(a+b-1)(n-1)}{b(n-1)-bk-2}|W|=\frac{(a+b-1)(n-1)}{b(n-1)-bk-2}(n-\delta(G)),
$$
which hints that
$$
\delta(G)\geq\frac{(a-1)n+b+bk+2}{a+b-1}.\eqno(2)
$$

In light of (2), $H=G-D$ and $|D|=k$, we possess
$$
\delta(H)\geq\delta(G)-k\geq\frac{(a-1)n+b-(a-1)k+2}{a+b-1}.
$$
We finish the proof of Claim 1. \hfill $\Box$

Note that $\beta=\min\{d_{H-S}(x):x\in T\}$. Then there exists $x_1\in T$ such that $d_{H-S}(x_1)=\beta$. By virtue of
Claim 1 and $\delta(H)\leq d_H(x_1)\leq d_{H-S}(x_1)+|S|=\beta+|S|$, we acquire
$$
|S|\geq\delta(H)-\beta\geq\frac{(a-1)n+b-(a-1)k+2}{a+b-1}-\beta.\eqno(3)
$$

The following proof is divided into four cases by virtue of the value of $\beta$.

\noindent{\bf Case 1.} \ $\beta=a$.

It follows from (1) and $\varepsilon(S)\leq|S|$ that
\begin{eqnarray*}
\varepsilon(S)-1&\geq&\theta_H(S,T)=b|S|+d_{H-S}(T)-a|T|\\
&\geq&b|S|+\beta|T|-a|T|=b|S|\geq|S|\geq\varepsilon(S),
\end{eqnarray*}
it is a contradiction.

\noindent{\bf Case 2.} \ $a-1\geq\beta\geq2$.

By virtue of (1), (3), $|S|+|T|+k\leq n$ and $\varepsilon(S)\leq2$, we acquire
\begin{eqnarray*}
1&\geq&\varepsilon(S)-1\geq\theta_H(S,T)=b|S|+d_{H-S}(T)-a|T|\\
&\geq&b|S|+\beta|T|-a|T|=b|S|-(a-\beta)|T|\\
&\geq&b|S|-(a-\beta)(n-k-|S|)\\
&=&(a+b-\beta)|S|-(a-\beta)(n-k)\\
&\geq&(a+b-\beta)\Big(\frac{(a-1)n+b-(a-1)k+2}{a+b-1}-\beta\Big)-(a-\beta)(n-k).
\end{eqnarray*}

Let $\varphi(\beta)=(a+b-\beta)(\frac{(a-1)n+b-(a-1)k+2}{a+b-1}-\beta)-(a-\beta)(n-k)$. Then we get
$$
\varphi(\beta)\leq1.\eqno(4)
$$

It follows from $a-1\geq\beta\geq2$ and $n\geq\frac{(a+b-2)(2a+b-3)+2}{b}+\frac{bk}{b-1}\geq\frac{(a+b-2)(2a+b-3)+2}{b}+k$
that
\begin{eqnarray*}
\frac{d\varphi}{d\beta}&=&-\Big(\frac{(a-1)n+b-(a-1)k+2}{a+b-1}-\beta\Big)-(a+b-\beta)+(n-k)\\
&=&2\beta+\frac{b(n-k)-b-2}{a+b-1}-a-b\\
&\geq&4+\frac{(a+b-2)(2a+b-3)+2-b-2}{a+b-1}-a-b\\
&=&a-1+\frac{1}{a+b-1}>0.
\end{eqnarray*}
Hence, we easily see that $\varphi(\beta)$ attains its minimum value at $\beta=2$, that is,
$$
\varphi(\beta)\geq\varphi(2).\eqno(5)
$$

According to (4), (5) and $n\geq\frac{(a+b-2)(2a+b-3)+2}{b}+\frac{bk}{b-1}\geq\frac{(a+b-2)(2a+b-3)+2}{b}+k$, we deduce
\begin{eqnarray*}
1&\geq&\varphi(\beta)\geq\varphi(2)=(a+b-2)\Big(\frac{(a-1)n+b-(a-1)k+2}{a+b-1}-2\Big)-(a-2)(n-k)\\
&=&\frac{b(n-k)-(a+b-2)(2a+b-4)}{a+b-1}\\
&\geq&\frac{(a+b-2)(2a+b-3)+2-(a+b-2)(2a+b-4)}{a+b-1}\\
&=&\frac{a+b}{a+b-1}>1,
\end{eqnarray*}
a contradiction.

\noindent{\bf Case 3.} \ $\beta=1$.

\noindent{\bf Subcase 3.1.} \ $|T|>\big\lfloor\frac{(b(n-1)-bk)n-2(n-1)}{(a+b-1)(n-1)}\big\rfloor$.

In terms of the integrity of $|T|$, we have
$$
|T|\geq\bigg\lfloor\frac{(b(n-1)-bk)n-2(n-1)}{(a+b-1)(n-1)}\bigg\rfloor+1.\eqno(6)
$$

Note that $d_{H-S}(x_1)=\beta=1$ and $x_1\in T$. Then we easily see that
$$
x_1\notin N_G(T\setminus N_G(x_1)).\eqno(7)
$$

It follows from (6), $T\cap(S\cup D)=\emptyset$ and $H=G-D$ that
\begin{eqnarray*}
|T\setminus N_G(x_1)|&=&|T\setminus N_{G-D}(x_1)|=|T\setminus N_H(x_1)|=|T\setminus N_{H-S}(x_1)|\\
&\geq&|T|-|N_{H-S}(x_1)|=|T|-d_{H-S}(x_1)=|T|-1\\
&\geq&\bigg\lfloor\frac{(b(n-1)-bk)n-2(n-1)}{(a+b-1)(n-1)}\bigg\rfloor.
\end{eqnarray*}
Combining this with the assumption of Theorem 2, we gain
$$
N_G(T\setminus N_G(x_1))=V(G),
$$
which contradicts (7).

\noindent{\bf Subcase 3.2.} \ $|T|\leq\big\lfloor\frac{(b(n-1)-bk)n-2(n-1)}{(a+b-1)(n-1)}\big\rfloor$.

\noindent{\it Claim 2.} \ $|T|\leq\frac{b(n-1)-bk-2}{a+b-1}$.

\noindent{\it Proof.} \ Let $|T|>\frac{b(n-1)-bk-2}{a+b-1}$. According to (3) and $\beta=1$, we gain
$$
|S|+|T|+k>\frac{(a-1)n+b-(a-1)k+2}{a+b-1}-1+\frac{b(n-1)-bk-2}{a+b-1}+k=n-1.
$$
Combining this with $|S|+|T|+k\leq n$ and the integrity of $|S|+|T|+k$, we have
$$
|S|+|T|+k=n.\eqno(8)
$$

From (8), $\varepsilon(S)\leq2$ and $|T|\leq\big\lfloor\frac{(b(n-1)-bk)n-2(n-1)}{(a+b-1)(n-1)}\big\rfloor\leq\frac{(b(n-1)-bk)n-2(n-1)}{(a+b-1)(n-1)}$,
we refer
\begin{eqnarray*}
\theta_H(S,T)&=&b|S|+d_{H-S}(T)-a|T|\\
&\geq&b|S|+|T|-a|T|\\
&=&b(n-k-|T|)-(a-1)|T|\\
&=&b(n-k)-(a+b-1)|T|\\
&\geq&b(n-k)-(a+b-1)\cdot\frac{(b(n-1)-bk)n-2(n-1)}{(a+b-1)(n-1)}\\
&=&\frac{nbk}{n-1}-bk+2\geq2\geq\varepsilon(S),
\end{eqnarray*}
this contradicts (1). This finishes the proof of Claim 2. \hfill $\Box$

Let $\lambda=|\{x:x\in T, d_{H-S}(x)=1\}|$. Evidently, $\lambda\geq1$ and $|T|\geq\lambda$. Using (3), Claim 2, $b\geq a\geq2$,
$\beta=1$ and $\varepsilon(S)\leq2$, we acquire that
\begin{eqnarray*}
\theta_H(S,T)&=&b|S|+d_{H-S}(T)-a|T|\\
&\geq&b|S|+2|T|-\lambda-a|T|\\
&=&b|S|-(a-2)|T|-\lambda\\
&\geq&b\bigg(\frac{(a-1)n+b-(a-1)k+2}{a+b-1}-1\bigg)-(a-2)\cdot\frac{b(n-1)-bk-2}{a+b-1}-\lambda\\
&=&2+\frac{b(n-1)-bk-2}{a+b-1}-\lambda\geq2+|T|-\lambda\geq2\geq\varepsilon(S),
\end{eqnarray*}
this conflicts with (1).

\noindent{\bf Case 4.} \ $\beta=0$.

Let $d=|\{x:x\in T, d_{H-S}(x)=0\}|$. Distinctly, $d\geq1$. Setting $Z=V(H)\setminus S$. Thus, $N_G(Z)\neq V(G)$ since $\beta=0$.
By the assumption of Theorem 2, we get
$$
n-d\geq|N_G(Z)|\geq\frac{(a+b-1)(n-1)}{b(n-1)-bk-2}|Z|=\frac{(a+b-1)(n-1)}{b(n-1)-bk-2}(n-k-|S|),
$$
which implying
$$
|S|\geq n-k-\frac{(n-d)(b(n-1)-bk-2)}{(a+b-1)(n-1)}. \eqno(9)
$$

By virtue of $n\geq\frac{(a+b-2)(2a+b-3)+2}{b}+\frac{bk}{b-1}$, we easily verify that
$$
\frac{b(n-1)-bk-2}{n-1}>1. \eqno(10)
$$

Using (9), (10), $|S|+|T|+k\leq n$, $b\geq a\geq2$ and $\varepsilon(S)\leq2$, we deduce
\begin{eqnarray*}
\theta_H(S,T)&=&b|S|+d_{H-S}(T)-a|T|\\
&\geq&b|S|+|T|-d-a|T|\\
&=&b|S|-(a-1)|T|-d\\
&\geq&b|S|-(a-1)(n-k-|S|)-d\\
&=&(a+b-1)|S|-(a-1)(n-k)-d\\
&\geq&(a+b-1)\bigg(n-k-\frac{(n-d)(b(n-1)-bk-2)}{(a+b-1)(n-1)}\bigg)-(a-1)(n-k)-d\\
&=&b(n-k)-\frac{(n-d)(b(n-1)-bk-2)}{n-1}-d\\
&\geq&b(n-k)-\frac{(n-1)(b(n-1)-bk-2)}{n-1}-1\\
&=&b+1>2\geq\varepsilon(S),
\end{eqnarray*}
which conflicts with (1). Theorem 2 is justified. \hfill $\Box$

\medskip

\section{Remark}

Now, we show that the condition on neighborhood in Theorem 2 is sharp, that is, we cannot replace it by $N_G(X)=V(G)$ or
$$
|N_G(X)|\geq\frac{(a+b-1)(n-1)}{b(n-1)-bk-2}|X|
$$
for all $X\subseteq V(G)$.

Let $a,b,k$ and $t$ be nonnegative integers such that $b\geq a\geq2$, $t$ is sufficiently large, and $\frac{t+1}{2}$ and
$\frac{(a-1)t+2}{b}$ are two integers. We construct a graph $G=K_{\frac{(a-1)t+2}{b}+k}\vee(\frac{t+1}{2}K_2)$. Setting
$n=|V(G)|=\frac{(a-1)t+2}{b}+k+t+1$, $A=V(K_{\frac{(a-1)t+2}{b}+k})$, $B=V(\frac{t+1}{2}K_2)$ and $D\subseteq A$ with
$|D|=k$. Next, we verify that $N_G(X)=V(G)$ or
$$
|N_G(X)|\geq\frac{(a+b-1)(n-1)}{b(n-1)-bk-2}|X|
$$
holds for all $X\subseteq V(G)$.

We easily see that $N_G(X)=V(G)$ if $|X\cap A|\geq2$, or $|X\cap A|=1$ and $|X\cap B|\geq1$. Of course, if
$|X|=1$ and $X\subseteq A$, then we easily claim
$$
|N_G(X)|=|V(G)|-1=n-1>\frac{(a+b-1)(n-1)}{b(n-1)-bk-2}|X|.
$$
Hence, we may assume that $X\subseteq B$. In this case, we possess
$$
|N_G(X)|=|A|+|X|=\frac{(a-1)t+2}{b}+k+|X|.
$$
Therefore,
$$
|N_G(X)|\geq\frac{(a+b-1)(n-1)}{b(n-1)-bk-2}|X|
$$
holds if and only if
$$
\frac{(a-1)t+2}{b}+k+|X|\geq\frac{(a+b-1)(n-1)}{b(n-1)-bk-2}|X|.
$$
This inequality above is equivalent to $|X|\leq t$. Thus if $X\neq B$ and $X\subset B$, then we derive
$$
|N_G(X)|\geq\frac{(a+b-1)(n-1)}{b(n-1)-bk-2}|X|.
$$
If $X=B$, then we easily see that $N_G(X)=V(G)$. Consequently, $N_G(X)=V(G)$ or
$$
|N_G(X)|\geq\frac{(a+b-1)(n-1)}{b(n-1)-bk-2}|X|
$$
holds for all $X\subseteq V(G)$.

Finally, we demonstrate that $G$ is not fractional $(a,b,k)$-critical covered. Let $H=G-D$. For above $A\setminus D$ and $B$,
we admit $|A\setminus D|=\frac{(a-1)t+2}{b}$, $|B|=t+1$, $d_{H-(A\setminus D)}(B)=t+1$ and $\varepsilon(A\setminus D)=2$. Thus,
we acquire
\begin{eqnarray*}
\theta_H(A\setminus D,B)&=&b|A\setminus D|+d_{H-(A\setminus D)}(B)-a|B|\\
&=&(a-1)t+2+t+1-a(t+1)=3-a\leq1<2=\varepsilon(A\setminus D).
\end{eqnarray*}
By Theorem 3, $H$ is not fractional $[a,b]$-covered, and so, $G$ is not fractional $(a,b,k)$-critical covered.

\medskip


\end{document}